\documentclass[12pt]{amsart}
\calclayout
\newcommand{\SD}{\mathbb{SD}}
\newcommand{\SPD}{\mathbb{SPD}}

\newcommand{\tauiso}{\tau^{iso}}
\newcommand{\di}{\displaystyle}
\calclayout 
\usepackage{amssymb}
\usepackage[all]{xy}
\usepackage{amsthm}
\usepackage{amsfonts}
\usepackage{amsmath}
\usepackage{amscd}
\newtheoremstyle{andy}{0pt}{0pt}{\itshape}{}{\bf\sf}{.}{4pt}{}
\newtheorem{claim}{Claim}

\newtheorem{lemma}[claim]{Lemma}

\newtheorem{prop}[claim]{Proposition}
\newtheorem{defin}[claim]{Definition}
\newtheorem{corollary}[claim]{Corollary}
\newtheorem{sign}[claim]{Sign Convention}

\newtheorem{example}[claim]{Example}
\newtheorem{assumption}[claim]{Assumption}

\newcommand{\bA}{\mathbb{A}}

\newcommand{\C}{\mathbf{C}}

\newcommand{\Z}{\mathbf{Z}}
\newcommand{\Q}{\mathbf{Q}}

\newcommand{\infinity}{\infty}
\newcommand{\CC}{\mathcal{C}}

\newcommand{\rR}{\mathrm{rank}_R}
\newcommand{\taunew}{\tau^{NEW}}

\newcommand{\Image}{\mathrm{Im}}
\newcommand{\newt}{\taunew}
\newcommand{\taui}{\tau^{iso}}
\newcommand{\Kiso}{K_1^{iso}}
\newcommand{\CP}{\mathbf{CP}}
\newcommand{\taunewiso}{\taunew_{iso}}

\begin{document}

\title {Absolute Whitehead Torsion}
\author{Andrew Korzeniewski}
\address{School of Mathematics
 \newline\indent
University of Edinburgh
 \newline\indent
Edinburgh EH9 3JZ, United Kingdom}
\email{A.J.Korzeniewski@sms.ed.ac.uk}

\begin{abstract}
We refine the Whitehead torsion of a chain equivalence of finite
chain complexes in an additive category $\bA$ from an element of
$\widetilde{K}^{iso}_1(\bA)$ to an element of the absolute group
$K_1^{iso}(\bA)$.  We apply this invariant to symmetric Poincar\'e
complexes and identify it in terms of more traditional invariants.
In the companion paper \cite{hkr} (joint with Ian Hambleton and
Andrew Ranicki) this new invariant is applied to obtain the
multiplicativity of the signature of fibre bundles mod 4.
\end{abstract}

\maketitle

\section*{Introduction}
The Whitehead torsion of a homotopy equivalence $f:X \to Y$ of
finite $CW$ complexes is an element of the Whitehead group of
$\pi=\pi_1(X)=\pi_1(Y)$
$$\tau(f)~=~\tau(\widetilde{f}:C(\widetilde{X}) \to C(\widetilde{Y}))
 \in Wh(\pi)~=~K_1({\mathbb Z}[\pi])/\{\pm \pi\}~,$$
with $\widetilde{f}$ the induced chain equivalence of based f.g. free
cellular ${\mathbb Z}[\pi]$-module chain complexes.
The Whitehead torsion of a finite $n$-dimensional Poincar\'e complex $X$ is
$$\tau(X)~=~\tau([X] \cap-:C(\widetilde{X})^{n-*} \to C(\widetilde{X}))
\in Wh(\pi)~.$$
In this paper we extend the methods of \cite{alg} to consider absolute
Whitehead torsion invariants for homotopy equivalences of certain finite
$CW$ complexes and finite Poincar\'e complexes, which take values in
$K_1({\mathbb Z}[\pi])$ rather than $Wh(\pi)$. We shall also be
extending the round $L$-theory of \cite{rltheory}, which is the algebraic
$L$-theory with absolute Whitehead torsion decorations.
In the paper \cite{hkr} absolute Whitehead torsion in both
algebraic
$K$- and $L$-theory will be applied to investigate the signatures of
fibre bundles.

The absolute torsion of a finite contractible chain complex of
finitely generated based $R$-modules
 $C$ is defined by
$$\tau(C)~=~\tau(d+\Gamma:C_{odd} \to C_{even}) \in K_1(R)~.$$
for a chain contraction $\Gamma$; it is independent of the choice of
$\Gamma$.
The algebraic mapping cone of a chain equivalence of finite chain
complexes of finitely generated based $R$-modules $f:C \to D$ is a
contractible chain complex $\CC(f)$.
The naive absolute torsion $\tau(C(f)) \in K_1(R)$ only has good additive
and composition formulae modulo ${\rm im}(K_1({\mathbb Z}) \to K_1(R))$.
Likewise, the naive definition of the torsion of an $n$-dimensional
symmetric Poincar\'e complex $(C,\phi)$ (see \cite{origandrew}) with
$C$ a based f.g. free $R$-module chain complex
$$\tau(C,\phi)~=~\tau(\CC(\phi_0:C^{n-*} \to C))$$
only has good cobordism and additivity properties in
$\widetilde{K}_1(R)$.  The Tate $\Z_2$-cohomology class
\[\tau(C,\phi) \in \widehat{H}^n(\Z_2;K_1(R))\]
may not be defined, and even if defined may not be a cobordism
invariant.

In \cite{alg} Andrew Ranicki developed a theory of absolute torsion
for chain equivalences of round chain complexes, that is chain
complexes $C$ satisfying $\chi(C)=0$.  This absolute torsion has a
good composition formula but it is not additive, and for round
Poincar\'e complexes $\tau(C,\phi_0)$ is not a cobordism invariant
(contrary to the assertions of \cite[7.21, 7.22]{andrew:additive}).

There are two main aims of this paper, firstly to develop a more
satisfactory definition of the absolute torsion of a chain equivalence
with good additive and composition formulae and secondly to define an
absolute torsion invariant of Poincar\'e complexes which behaves
predictably under cobordism.  Chapter \ref{section:abs} is devoted to
the first of these aims.  Following \cite{alg} we work in the more
general context of an additive category $\bA$.  The chief novelty here is the
introduction
of a \textit{signed} chain complex; this is a pair $(C,\eta_C)$ where
$C$ is a finite chain complex and $\eta_C$ is a ``sign'' term living
in $K_1(\bA)$, which will be made precise in chapter \ref{section:abs}.
We give definitions for the sum and suspension of two signed chain
complexes, and we define the absolute torsion of a chain equivalence
of signed chain complexes.
\[\taunew(f:C \to D) \in K_1^{iso}(\bA)\]
 This gives us a definition of absolute
torsion with good additive and composition formulae at the cost of
making the definition more complicated by adding sign terms to the chain
complexes.  This definition is similar to the one given in \cite{alg},
indeed if the chain complexes $C$ and $D$ are round and
$\eta_C=\eta_D=0$ then the definition of the absolute torsion of a
chain equivalence $f:C \to D$ is precisely that given in \cite{alg}.
When working over a ring $R$ the absolute torsion defined here
reduces to the usual torsion in $\widetilde{K_1}(R)$.

In chapter \ref{section:dual} we work over a category with involution
and define the dual of a signed chain complex.  We can then
define in chapter \ref{section:poincare} the absolute
torsion of a symmetric Poincar\'e complex to be the absolute torsion
of the chain equivalence $\phi_0:C^{n-*}\to C$.  This new invariant is
shown to be additive and to have good behaviour under \textit{round}
algebraic cobordism.  Although we have to choose a sign $\eta_C$ in
order to define the absolute torsion of $\phi_0$, we show that the
absolute torsion is independent of this choice.  In chapter
\ref{section:product} we state a product formula for the absolute torsion
and prove it for group rings.

It should be noted that the definitions in ``Round
$L$-Theory'' \cite{rltheory} used the absolute torsion invariant of
\cite{alg}, \cite{andrew:additive}, which is not a cobordism
invariant.  In chapter
\ref{section:rltheory} we show that the absolute torsion
defined here may be used as the ``correct'' definition; in this case
the statements in \cite{rltheory} are correct.

In chapter \ref{section:manifolds} we investigate the absolute
torsion of manifolds.  This invariant is only defined when we pass to
the reduced group $\widehat{H}^n(\Z_2;K_1(\Z[\pi_1M]))$ and
we provide some examples.  In chapter \ref{section:sign} the ``sign''
term of the absolute torsion of a manifold is identified with more
traditional invariants of a manifold such as the signature, the Euler
characteristic and the semi-characteristic.

The forthcoming paper \cite{hkr} will make extensive use of the
invariants and techniques developed here.  Chapters
\ref{section:derived} and \ref{section:spd} develop the notion of the
\textit{signed derived category} which will be required by \cite{hkr}.

I would like to thank my supervisor Andrew Ranicki for his help and
encouragement during the writing of this paper.  I would also like to
thank Ian Hambleton for many useful conversations and for carefully
checking the computations.

\tableofcontents

\section{Absolute torsion of contractible complexes and chain equivalences.}
\label{section:abs}

In this section we introduce the absolute torsion of contractible
complexes and chain equivalences and derive their basic properties.
This closely follows \cite{alg} but without the assumption that the
complexes are round (\(\chi(C) =0 \in K_0(\bA)\)); we also develop the
theory in the context of signed chain complexes which we will define
in this section.

Let $\bA$ be an additive category.  Following \cite{alg} we define:

\begin{defin}
\begin{enumerate}
\item The \emph{class group} $K_0(\bA)$
has one generator $[M]$ for each object in $\bA$ and relations:
\begin{enumerate}
\item $[M]=[M']$ if $M$ is isomorphic to $[M']$.
\item $[M\oplus N]=[M]+[N]$ for objects $M,N$ in $\bA$.
\end{enumerate}
\item The \emph{isomorphism torsion group} $K_1^{iso}(\bA)$ has one
generator $\tau^{iso}(f)$ for each isomorphism $f:M \to N$ in $\bA$,
and relations:
\begin{enumerate}
\item $\tau^{iso}(gf)=\taui(f)+\taui(g)$
 for isomorphisms
$f:M\to N$, $g:N\to P$
\item $\taui(f \oplus f') = \taui(f) + \taui(f')$
 for isomorphisms
$f:M \to N$, $f':M'\to N'$
\end{enumerate}
\item The \emph{automorphism torsion group} $K_1(\bA)$ has one
generator $\tau(f)$ for each automorphism $f:M \to N$ in $\bA$,
and relations:
\begin{enumerate}
\item $\tau(gf)=\tau(f)+\tau(g)$ for automorphisms $f:M \to N$,
$g:N \to P$.
\item $\tau(f \oplus f') = \tau(f) + \tau(f')$ for automorphisms
$f:M \to N$, $f':M' \to N'$.
\end{enumerate}
\end{enumerate}
\end{defin}

\subsection{Sign terms}
The traditional torsion invariants are considered to lie in
\(\widetilde{K}_1^{iso}(\bA)\), a particular quotient of
$K_1^{iso}(\bA)$ (defined below) in which the torsion of maps such as
\(\left(\begin{array}{cc} 0 & 1 \\ 1 & 0\end{array}\right):C \oplus D
\to D \oplus C\) are
trivial.
In absolute torsion we must consider such
rearrangement maps;
to this end we recall from \cite{alg} the following notation:
\begin{defin}
Let \(C,D\) be free, finitely generated chain
complexes in \(\bA\).
\begin{enumerate}
\item The \emph{suspension} of \(C\) is the chain complex \(SC\) such
that \(SC_r=C_{r-1}\) and \(SC_0=0\)
\item The \emph{sign} of two objects \(X,Y \in \bA\) is
the element
\[\epsilon(X,Y) := \taui\left(\left(\begin{array}{cc}0 & 1_Y \\ 1_X &
0\end{array}\right)\colon X\oplus Y \rightarrow Y \oplus X\right) \in
\Kiso(\bA)\]
The sign only depends on the stable isomorphism classes of $M$ and $N$
and satisfies:
\begin{enumerate}
\item $\epsilon(M \oplus M',N) = \epsilon(M,N)+\epsilon(M',N)$
\item $\epsilon(M,N) = -\epsilon(N,M)$
\item $\epsilon(M,M) = \taui(-1:M \to M)$
\end{enumerate}
We may extend $\epsilon$ to a morphism of abelian groups:
\[\epsilon\colon K_0(\bA)\otimes K_0(\bA)\to \Kiso(\bA);([M],[N]) \mapsto
\epsilon(M,N)\]
\item The \emph{reduced isomorphism torsion group}
$\widetilde{K}_1^{iso}(\bA)$ is the quotient:
\[\widetilde{K}_1^{iso}(\bA):= K^{iso}_1(\bA)/\Image(\epsilon\colon
K_0(\bA) \otimes K_0(\bA) \to K_1^{iso}(\bA))\]
\item The \emph{intertwining} of \(C\) and \(D\) is the element
defined by:
\[\beta(C,D):=\sum_{i>j}(\epsilon(C_{2i},D_{2j}) -
\epsilon(C_{2i+1},D_{2j+1})) \in\Kiso(\bA)\]
\end{enumerate}
\end{defin}

\begin{example}
The reader may find it useful to keep the following example in mind,
as it is the most frequently occurring context.

Let $R$ be an associative ring with 1 such that $\rR(M)$ is
well-defined for f.g
free modules $M$.
We define $\bA(R)$ to be the category of based f.g. $R$-modules.
In this case the map \(K_0(\bA(R)) \to \Z\) given by $M\mapsto\dim M$ is
an isomorphism.
We have a forgetful functor:
\[\Kiso(\bA(R)) \to K_1(R)\;;\;\taui(f) \mapsto \tau(f)\]
mapping elements of $\Kiso(\bA(R))$ to the more familiar $K_1(R)$ in the
obvious way.  In particular
\[\Image (\epsilon\colon K_0(\bA(R)) \otimes K_0(\bA(R)) \to K_1(R)) =
\{\tau(\pm 1)\} =\Image(K_1(\Z) \to K_1(R))\]
justifying the terminology of a ``sign'' term;
the map is given explicitly for modules $M$ and $N$ by:
\[\epsilon(M,N) = \rR(M)\rR(N)\tau(-1)\]

\end{example}

We will make use of the notation:
\[C_{even} = C_0 \oplus C_2 \oplus C_4 \oplus ...\]
\[C_{odd} = C_1 \oplus C_3 \oplus C_5 \oplus ...\]
and as usual we define the Euler characteristic $\chi(C)$ as:
\[\chi(C)=[C_{even}]-[C_{odd}] \in K_0(\bA)\]

We also recall from \cite{alg} proposition 3.4 the following
relationships between the ``sign'' terms:

\begin{lemma}
\label{lemma:beta}
Let \(C,C',D,D'\) be finite chain complexes over \(\bA\). Then
\begin{enumerate}
\item\(\beta(C,D) = \taui((C\oplus D)_{even} \to C_{even}\oplus
D_{even}) \\- \taui((C\oplus D)_{odd} \to C_{odd}\oplus
D_{odd})\)
\item \(\beta(C \oplus C',D) = \beta(C,D)+\beta(C',D)\)
\item \(\beta(C,D \oplus D') = \beta(C,D)+\beta(C,D')\)
\item \(\beta(C,D) - \beta(D,C) +  \sum (-)^r\epsilon(C_r,D_r) =
\epsilon(C_{even},D_{even}) - \epsilon(C_{odd},D_{odd})\)
\item \(\beta(SC,SD) = -\beta(C,D)\)
\item \(\beta(SC,C) = \epsilon(C_{odd},C_{even})\)
\end{enumerate}
\end{lemma}

\subsection{Signed chain complexes.}

In order to make to formulae in this paper more concise we introduce
the concept of a \emph{signed} chain complex; this is a chain complex
with an associated element in \(\Image(\epsilon\colon K_0(\bA)\otimes
K_0(\bA)\to \Kiso(\bA))\) which we
refer to as the \emph{sign} of the complex.  We use this
element in the definition of the absolute torsion invariants.

\begin{defin}
\begin{enumerate}
\item A \emph{signed} chain complex is a pair \((C,\eta_C)\) where
\(C\) is a finite chain complex in $\bA$ and \(\eta_C\) an
element of
\[\Image(\epsilon\colon K_0(\bA)\otimes
K_0(\bA)\to \Kiso(\bA))\]
  We will usually suppress
mention of \(\eta_C\) denoting such complexes as \(C\).
\item Given a signed chain complex  \((C,\eta_C)\) we give the
suspension of \(C\), \(SC\) the sign
\[\eta_{SC} = -\eta_C\]
\item We define the \emph{sum} signed chain complex of two signed
chain complexes \((C,\eta_C)\), \((D,\eta_D)\) as \((C\oplus
D,\eta_{C\oplus D})\) where \(C\oplus D\) is the usual based sum of
two chain complexes and \(\eta_{C\oplus D}\) defined by:
\[\eta_{C\oplus
D}=\eta_C+\eta_D-\beta(C,D)+\epsilon(C_{odd},\chi(D))\]
(it is easily shown that \(\eta_{(C\oplus D) \oplus E} = \eta_{C
\oplus (D \oplus E)}\))
\end{enumerate}
\end{defin}

\subsection{The absolute torsion of isomorphisms}
We now define the absolute torsion of a collection of isomorphisms
$\{f_r:C_r \to D_r\}$ between two signed chain complexes.  Note that the
map $f$ need not be a chain isomorphism (i.e. $fd_C=d_Df$ need not
hold).  In the case where $f$ is a chain isomorphism the torsion
invariant defined here will coincide with the definition of the
absolute torsion of chain equivalence given later.
\begin{defin}
\label{defin:iso}
The \emph{absolute torsion} of a collection of  isomorphisms
$\{f_r:C_r \to D_r\}$  between the chain groups of signed chain
complexes $C$ and $D$ is defined as:
\[\taunewiso(f) = \sum_{r=0}^{\infinity}(-)^r\taui(f_r:C_r\to D_r) -
\eta_C + \eta_D \in \Kiso(\bA)\]
\end{defin}
\begin{lemma}
\label{lemma:isomorphism}
We have the following properties of the absolute torsion of
isomorphisms:
\begin{enumerate}
\item The absolute torsion of isomorphisms is logarithmic, that is for
isomorphisms \(f:C \to D\) and \(f:D \to E\).
\[\taunewiso(gf) = \taunewiso(f) + \taunewiso(g)\]
\item The absolute torsion of isomorphisms is additive, that is for
isomorphisms \(f: C \to D\) and \(f':C'\to D'\)
\[\taunewiso(f\oplus g) = \taunewiso(f) + \taunewiso(g)\]
\item The absolute torsion of the rearrangement isomorphism:
\[ C \oplus D \xrightarrow{\left(\begin{array}{cc}0 & 1 \\ 1 &
0\end{array} \right)} D \oplus C\]
is \(\epsilon(\chi(C),\chi(D)) \in \Kiso(\bA)\).
\item The absolute torsion of the isomorphism:
\[ S(C \oplus D) \xrightarrow{\left(\begin{array}{cc}1 & 0 \\ 0 &
1\end{array} \right)} SC \oplus SD\]
is \(\epsilon(\chi(D),\chi(C))\in \Kiso(\bA)\).
\end{enumerate}
\end{lemma}
\begin{proof}
Parts 1 and 2 follow straight from the definitions.  For part 3
\begin{eqnarray*}
\taunewiso(C\oplus D \to D\oplus C) & = &
\sum_{r=0}^\infinity(-)^r\epsilon(C_r,D_r)-\eta_{C\oplus D} + \eta_{D
\oplus C} \\
& = & \sum_{r=0}^\infinity(-)^r\epsilon(C_r,D_r)+\beta(C,D) - \beta(D,C) \\
& & -
\epsilon(C_{odd},\chi(D)) + \epsilon(D_{odd},\chi(C)) \\
& = & \epsilon(C_{even},D_{even}) - \epsilon(C_{odd},D_{odd}) \\
& & -
\epsilon(C_{odd},\chi(D)) + \epsilon(D_{odd},\chi(C)) \\
& = & \epsilon(\chi(C),\chi(D))
\end{eqnarray*}
For part 4:
\begin{eqnarray*}
\taunewiso(S(C\oplus D) \to SC \oplus SD) & = & \eta_{SC\oplus SD} -
\eta_{S(C\oplus D)} \\
& = & -\beta(SC,SD) + \epsilon(C_{even},\chi(SD)) \\
& & - \beta(C,D) +
\epsilon(C_{odd},\chi(D)) \\
& = & \epsilon(\chi(D),\chi(C))
\end{eqnarray*}
\end{proof}

\subsection{The absolute torsion  of contractible complexes and short
exact sequences.}
We recall from \cite{alg} the following:

Given a finite contractible chain complex over \(\bA\)
\[C\colon C_n \rightarrow ...\rightarrow
C_0\]
and a chain contraction \(\Gamma\colon C_r\rightarrow C_{r+1}\) we may
form the following isomorphism:
\small
\[d+\Gamma\;=\;\left(\begin{array}{cccc}d&0&0&...\\\Gamma&d&0&...\\
0&\Gamma&d&...\\\vdots&\vdots&\vdots&\end{array}\right)\colon
C_{odd}=C_1\oplus C_3 \oplus C_5 ... \rightarrow C_{even}=C_0 \oplus
C_2 \oplus C_4...\]
\normalsize
The element \(\taui(d+\Gamma) \in \Kiso(\bA)\) is independent of the choice of
\(\Gamma\) and is denoted \(\tau(C)\) (following \cite{alg} section 3).

We define the \emph{absolute torsion} of a contractible signed chain complex
\(C\) as
\[\taunew(C) = \tau(C) + \eta_C \in \Kiso(\bA)\]

Given a short exact sequence of signed chain complexes over $\bA$:
\[0 \to C \xrightarrow{i} C'' \xrightarrow{j} C' \to 0\]
we may find a sequence of splitting morphisms $\{k:C'_r \to C''_r|r
\geq 0 \}$ such that $jk=1:C'_r \to C'_r\;(r \geq 0)$ and each
$(i\;k):C_r\oplus C'_r \to C''_r\;(r \geq 0)$ is an isomorphism.  The
torsion of this collection of isomorphisms
\[\taunewiso((i\;k):C_r\oplus
C'_r \to C''_r)\]
 is independent of the choice of the $k_r$, so we may define
the \emph{absolute torsion of a short exact sequence} as:
\[\taunew(C,C'',C';i,j) = \taunewiso((i\;k):C_r\oplus
C'_r \to C''_r)\]

\begin{lemma}
\label{lemma:contract}
We have the following properties of the absolute torsion of
signed contractible complexes:
\begin{enumerate}
\item
\label{lemma:contract:shortexact}
Suppose we have a short exact sequence of contractible signed complexes:
\[0\to C \xrightarrow{i} C'' \xrightarrow{j} C'\to 0\]
Then
\[\taunew(C'') = \taunew(C)+\taunew(C')+\taunew(C,C'',C';i,j)\]
\item
\label{lemma:contract:sum}
Let \(C\), \(C'\) be contractible signed complexes.  Then:
\[\taunew(C \oplus C') = \taunew(C) + \taunew(C')\]
\end{enumerate}
\end{lemma}
\begin{proof}
\begin{enumerate}
\item From \cite{alg} proposition 3.3 we have that
\[\tau(C'') = \tau(C) + \tau(C') +
\sum_{r=0}^\infinity\taui((i\;k):C_r\oplus C'_r \to C'') +
\beta(C,C)\]
for some choice of splitting morphisms $\{k:C'_r \to C''_r|r
\geq 0 \}$.
By the definition of the absolute torsion of a short exact sequence and
the definition of the sum torsion (noting that contractible complexes
have \(\chi(C)=0 \in K_0(\bA)\)) we get:
\begin{eqnarray*}
\taunew(C,C'',C';i,j) &=& \sum_{r=0}^\infinity(-)^r\taui((i\;k):C_r\oplus C'_r \to
C'')\\&& + \beta(C,C') - \eta_C - \eta_{C'} + \eta_{C''}
\end{eqnarray*}
By comparing these two formulae and the definition of the absolute
torsion of a contractible signed complex, the result follows.
\item Apply the above to \(C''=C \oplus C'\).
\end{enumerate}
\end{proof}
\subsection{The absolute torsion of chain equivalences.}
\begin{sign}
We define the algebraic mapping cone of a chain map \(f:C\rightarrow
D\) as
follows:
\[d_{\CC (f)} = \left(\begin{array}{cc} d_D & (-)^{r+1}f \\
0 & d_C \end{array}\right) \;:\;\CC(f)_r = D_r \oplus C_{r-1}
\rightarrow \CC(f)_{r-1} = D_{r-1} \oplus C_{r-2}\]
We make \(\CC(f)\) into a signed complex by setting
\[\eta_{\CC(f)} = \eta_{D \oplus SC}\]
\end{sign}

\begin{lemma}
The absolute torsion of a chain isomorphism \(f:C \to D\) of signed chain
complexes satisfies:
\[\taunewiso(f) = \taunew(\CC(f))\]
\end{lemma}
\begin{proof}
In the case of an isomorphism we may choose the chain contraction for
\(\CC(f)\) to be:
\[\Gamma_{\CC(f)} = \left(\begin{array}{cc} 0 & 0 \\ (-)^{r}f^{-1} &
0\end{array}\right) : \CC(f)_r \to \CC(f)_{r+1}\]
We have a commutative diagram:
\small
\[
\xymatrix{
(D_1 \oplus C_0) \oplus (D_3 \oplus C_2) \oplus ...
\ar[rrr]^*{(d_{\CC(f)}+\Gamma_{\CC(f)})} \ar[ddd] & & &
D_0 \oplus (D_2 \oplus C_1) \oplus (D_4 \oplus C_3)\oplus ...
\ar[ddd] \\ \\ \\
C_0 \oplus D_1 \oplus C_2 \oplus D_3...
\ar[rrr]^*+{\left(\begin{array}{cccc}f & d_D & 0 & \hdots\\
0 & -f^{-1} & d_C & \hdots \\
0 & 0 & f & \hdots \\
\vdots & \vdots & \vdots & \end{array}\right)} & & &
D_0 \oplus C_1 \oplus D_2 \oplus C_3
}
\]
\normalsize
The torsion of the upper map is \(\taui(\CC(f))\), the torsion
of the lower isomorphism is \(\sum_{r=0}^\infinity(-)^r\taui(f_r:C_r\to
D_r)+\epsilon(C_{odd},C_{odd})\) and the difference between the
torsions of the
downward maps is \(\sum_{r=0}^\infinity(-)^r\epsilon(C_r,C_{r-1})\) (using
the fact that \(C_r \cong D_r\)).  Hence
\begin{eqnarray*}
\taunew(\CC(f)) & = & \taui(\CC(f)) + \eta_{\CC(f)} \\
& = & \sum_{r=0}^\infinity(-)^r\taui(f_r:Cr\to
D_r) -
\sum_{r=0}^\infinity(-)^r\epsilon(C_r,C_{r-1}) \\
& & - \beta(C,SC) +
\epsilon(C_{odd},\chi(SC))+\epsilon(C_{odd},C_{odd})-\eta_C+\eta_D \\
& = & \sum_{r=0}^\infinity(-)^r\taui(f_r:C_r\to
D_r)-\eta_C+\eta_D \\
& = & \taunewiso(f)
\end{eqnarray*}
(using the formulae of lemma \ref{lemma:beta})
\end{proof}

We can now give a definition of the absolute torsion of a chain
equivalence $f:C \to D$ which coincides with the previous definition
in the case when $f$ is a chain isomorphism.

\begin{defin}
We define the \emph{absolute torsion} of a chain equivalence of signed
chain complexes \(f:C \rightarrow D\) as:
\[\taunew(f) = \taunew(\CC(f)) \in \Kiso(\bA)\]
In the case where \(f\) is a chain isomorphism the above lemma shows
that this definition of the torsion agrees with that given in
definition \ref{defin:iso}.
\end{defin}
\begin{lemma}
\label{lemma:equiv}
The absolute torsion of a chain equivalence of chain complexes
with torsion \(f: C\rightarrow D\) is:
\[\taunew(f) = \tau(\CC(f)) - \beta(D,SC) - \epsilon(D_{odd},\chi(C))
+ \eta_D - \eta_C \in K_1(\bA)\]
(c.f. definition of torsion on pages 223 and 226 of \cite{alg}.  The
two definitions coincide if \(C\) and \(D\) are even and \(\eta_C=\eta_D\)).
\end{lemma}
\begin{proof}
Simply a matter of unravelling definitions.
\end{proof}

We have the following properties of the torsion of chain equivalences:
\begin{prop}
\label{prop:torsionprops}
\begin{enumerate}
\item Let \(f:C \rightarrow D\) and \(g:D \rightarrow E\) be chain
equivalences of signed chain complexes in \(\bA\), then
\[\taunew(gf) = \taunew(f) + \taunew(g) \in \Kiso(\bA)\]
\item Suppose \(f:C \rightarrow D\) is map of contractible signed
chain complexes.
Then
\[\taunew(f) = \taunew(D) - \taunew(C) \in \Kiso(\bA)
\]

\item The absolute torsion \(\taunew(f)\) is a chain homotopy invariant
of \(f\).
\item
\label{prop:torsionprops:shortexact}
Suppose we have a commutative diagram of chain maps as follows where
the rows are
exact and the vertical maps are chain equivalences:
\[
\xymatrix{
0 \ar[r] & A \ar[r]^i \ar[d]_{a} & B \ar[r]^j \ar[d]_b &
C \ar[r] \ar[d]_c & 0 \\
0 \ar[r] & A' \ar[r]^{i'} & B' \ar[r]^{j'} & C' \ar[r] & 0
}
\]
Then
\begin{eqnarray*}
\taunew(b) &=& \taunew(a) + \taunew(c) - \taunew(A,B,C;i,j) \\
&& + \taunew(A',B',C';i',j') \in \Kiso(\bA)
\end{eqnarray*}
\item The torsion of a sum \(f \oplus f':C\oplus C' \rightarrow D
\oplus D'\) is given by:
\[
\taunew(f \oplus f') = \taunew(f) + \taunew(f')\in \Kiso(\bA)
\]
\item
\label{prop:torsionprops:contractexact}
Suppose we have a short exact sequence
\[0 \to A \xrightarrow{f} B \xrightarrow{g} C\to 0\]
where $C$ is a contractible complex and $f$ is a chain equivalence.
Then
\[\taunew(f) = \taunew(A,B,C;f,g) + \taunew(C)\]
\end{enumerate}
\end{prop}
\begin{proof}
The proofs of these follow those in \cite{alg} propositions 4.2 and
4.4, modified where appropriate.
\begin{enumerate}
\item We denote by \(\Omega C\) the chain complex defined by:
\[d_{\Omega C}=d_C:\Omega C_r = C_{r+1}\to \Omega C_{r-1}=C_r\]
We define a chain map
\[h:\Omega \CC(g) \to \CC(f)\]
by
\[\left(\begin{array}{cc} 0 & -1 \\ 0 &
0\end{array}\right):\Omega\CC(g)_r=E_{r+1}\oplus D_r \to \CC(f)_r =
D_r \oplus C_{r-1}\]
The algebraic mapping cone \(\CC(h)\) fits into the following short
exact sequences:
\begin{equation}
\label{eqn:shortone}
0\to\CC(f)\xrightarrow{i}\CC(h)\xrightarrow{j}\CC(g)\to 0
\end{equation}
\begin{equation}
\label{eqn:shorttwo}
0 \to \CC(gf)\xrightarrow{i'}\CC(h)\xrightarrow{j'}\CC(-1_D:D \to
D)\to 0
\end{equation}
where
\[i=\left(\begin{array}{c}1 \\
0\end{array}\right)\colon\CC(f)_r\to\CC(h)_r=\CC(f)_r\oplus\CC(g)_r\]
\[j=\left(\begin{array}{cc}0 &
1\end{array}\right)\colon\CC(h)_r=\CC(f)_r\oplus\CC(g)_r\to\CC(g)_r\]
\[i'=\left(\begin{array}{cc}0&0\\0&1\\1&0\\0&f\end{array}\right)
\colon\CC(gf)_r=E_r\oplus
SC_r \to \CC(h)_r=D_r\oplus SC_r\oplus E_r\oplus SD_r\]
\[j'=\left(\begin{array}{cccc}1&0&0&0\\0&-f&0&1\end{array}\right)
\colon\CC(h)_r=D_r\oplus
SC_r\oplus E_r\oplus SD_r \to \CC(-1_D)_r = D_r\oplus SD_r\]

Applying lemma \ref{lemma:contract} part
\ref{lemma:contract:shortexact} to the first short exact sequence
(\ref{eqn:shortone}) we
have
\begin{equation}
\label{eqn:f_g_is_h}
\taunew(h) = \taunew(f)+\taunew(g)
\end{equation}
Notice that
\begin{eqnarray*}
\taunewiso((i'\;k')) & = &
\taunewiso(\left(\begin{array}{cccc}0&0&1&0\\0&1&0&0\\
1&0&0&0\\0&f&0&1\end{array}\right)\colon\CC(gf)_r\oplus\CC(-1_D)_r \\
& & = E_r\oplus SC_r
\oplus D_r \oplus SD_r \to \CC(h)_r=D_r\oplus SC_r \oplus E_r\oplus SD_r) \\
& = & \taunewiso(D\oplus SC\to SC \oplus D) \\
& & + \taunewiso(E\oplus SC\to SC\oplus E) \\
& & + \taunewiso(D \oplus E \to E \oplus D) \\
& = & \epsilon(\chi(D),\chi(D))
\end{eqnarray*}
(using the results of lemma \ref{lemma:isomorphism}, the fact that
\(\chi(C)=\chi(D)=\chi(E)\) and that \(f\) has no effect on the torsion).
We also see that \(\taunew(\CC(-1_D)) = \taunewiso(-1_D) =
\epsilon(\chi(D),\chi(D))\).
Applying these two expressions and lemma \ref{lemma:contract} part
\ref{lemma:contract:shortexact} to the second exact sequence
(\ref{eqn:shorttwo}) we see that
\[\taunew(gf) = \taunew(h)\]
 and comparison with (\ref{eqn:f_g_is_h}) yields the result.
\item By construction we have \(\CC(0 \xrightarrow{0}D)=D\) and hence
\[\taunew(0 \xrightarrow{0} D) = \taunew(D)\]
Applying this and the composition formula (part 1) to the composition
\[0 \xrightarrow{0} C \xrightarrow{f} D\]
yields the result.
\item A chain homotopy
\[g\colon f \simeq f'\colon C\to D\]
gives rise to an isomorphism
\[\left(\begin{array}{cc}1 & (-)^{r}g \\ 0 &
1\end{array}\right)\colon\CC(f)=D \oplus CS \to \CC(f')=D \oplus SC\]
which has trivial torsion.  Using part 2
\[0= \taunew(\CC(f')) - \taunew(\CC(f))\]
the result follows.

\item
We choose splitting morphisms $\{k:C_r \to B_r|r
\geq 0 \}$ and $\{k':C'_r \to B'_r|r
\geq 0 \}$.
We have the following short exact sequence of mapping cones:
\[0 \rightarrow \CC(a) \xrightarrow{\left(\begin{array}{cc} i' & 0 \\ 0 &
i\end{array}\right)} \CC(b) \xrightarrow{\left(\begin{array}{cc} j' & 0 \\ 0 &
j\end{array}\right)} \CC(c) \rightarrow 0\]

We note that

\begin{eqnarray*}
&&\taunew\left(\CC(a),\CC(b),\CC(c);\left(\begin{array}{l}i'
 \\i\end{array}\right);
 \left(\begin{array}{l}i'\\i\end{array}\right)\right) \\
 &=&\taunewiso\left(\left(\begin{array}{cccc} i'
 & 0 & k' & 0 \\ 0 & i & 0 &
k\end{array}\right):A' \oplus SA \oplus C' \oplus SC \rightarrow
B'\oplus B\right) \\
& = & \taunewiso\left(\left(\begin{array}{cccc} i'
 & k' & 0  & 0 \\ 0 & 0 & i &
k\end{array}\right):A' \oplus C' \oplus SA \oplus SC \rightarrow
B'\oplus B\right) \\
& &+\taunewiso(SA\oplus C' \to C' \oplus SA)  \\
& = & \taunewiso((i\;k):SA\oplus SC \to SB) + \taunewiso((i'\;k')) +
\epsilon(\chi(SA),\chi(C')) \\
& = & \taunewiso((i\;k):S(A\oplus C) \to SB) + \taunewiso((i'\;k')) \\
&&+\taunewiso(SA\oplus SC\to S(A\oplus C)) + \epsilon(\chi(C),\chi(A)) \\
& = & - \taunewiso((i\;k):A\oplus C \to B) + \taunewiso((i'\;k')) \\
&& - \epsilon(\chi(C),\chi(A)) + \epsilon(\chi(C),\chi(A))\\
& = & - \taunewiso((i\;k):A\oplus C \to B) + \taunewiso((i'\;k')) \\
&=&\taunew(A',B',C';i',j') - \taunew(A,B,C;i,j)
\end{eqnarray*}

The result now follows from applying lemma \ref{lemma:contract} part
\ref{lemma:contract:shortexact} to the short exact sequence above.

\item Applying the result for a commutative diagram of short exact
sequences (part \ref{prop:torsionprops:shortexact}) with \(a=f:C \to
D\), $c=f':C' \to D'$ and $b=f \oplus f': C \oplus C' \to D \oplus D'$
yields the result.
\item We have a commutative diagram with short exact rows:
\[
\xymatrix{
0 \ar[r] & A \ar[r]^f \ar[d]_{f} & B \ar[r]^g \ar[d]_1 &
C \ar[r] \ar[d] & 0 \\
0 \ar[r] & B \ar[r]^{1} & B \ar[r] & 0 \ar[r] & 0
}
\]
The result follows by applying part \ref{prop:torsionprops:shortexact}
to the above diagram.

\end{enumerate}
\end{proof}

\subsection{Applications to topology and examples of use.}
\label{cp2_self_map}
Let $X$ be a connected finite $CW$-complex.  We may form the cellular chain
complex of the universal cover of $X$ as a complex $C(\tilde X)$ over
the fundamental group ring $\Z[\pi_1X]$; we may further make $C(\tilde
X)$ into a signed complex with an arbitrary choice of $\eta_{C(\tilde
X)}$. For a cellular homotopy
equivalence $f:X \to X$ we have an associated chain equivalence
$f_*:C(\tilde X) \to C(\tilde X)$; we can make $C(\tilde X)$ into a
signed chain complex by choosing some $\eta_{C(\tilde X)}$ and define
the torsion of $f$ to be
\[\taunew(f) := \taunew(f_*:C(\tilde X) \to C(\tilde X)) \in K_1(\Z[\pi_1X])\]
this is independent of the choice of $\eta_{C(\tilde X)}$.  We now
give some examples:

\begin{enumerate}
\item The torsion of the identity map of any connected $CW$-complex is
trivial.
\item Let $X = \CP^2$; we choose homogeneous coordinates $(x:y:z)$ and
we give $X$ a $CW$-structure as follows:

\begin{tabular}{ll}
0-cell & $(1:0:0)$ \\
2-cell & $(z_1:1:0)$ \\
4-cell & $(z_1:z_2:1)$ \\
\end{tabular}

Let $f:\CP^2 \to \CP^2$ be the cellular self-homeomorphism given by
complex conjugation in all three coordinates, that is:
\[f: (x:y:z) \mapsto (\bar x: \bar y:\bar z)\]
This map preserves the orientation of the 0-cell and 4-cell, and it
reverses the orientation of the 2-cell.  Hence $\taunew(f)=\tau(-1)$.
In corollary \ref{cor:selfmaps}
we show that for any orientation preserving self-homeomorphism
$g$ of a simply connected manifold of dimension $4k+2$, that
$\taunew(g)=0$.  This
example shows that for self-homeomorphism $f$ of a $4k$-dimensional
manifold it is possible for $\taunew(f) \neq 0$
\end{enumerate}

\section{The signed derived category.}
\label{section:derived}

The forthcoming paper \cite{hkr} will require the use of the signed
derived category $\SD(\bA)$.  In this section we define $\SD(\bA)$ and
prove some basic properties.

\begin{defin}
The \textit{signed derived category} $\SD(\bA)$ is the category with
objects signed chain complexes in $\bA$ and morphisms chain
homotopy classes of chain maps between such complexes.
\end{defin}

\begin{prop}\label{derivedprop}
{\rm (i)} The Euler characteristic defines a surjection
$$\chi~:~
K_0(\SD(\bA)) \to K_0(\bA)~;~[C,\eta_C] \mapsto
\chi(C)~=~\sum\limits^{\infty}_{r=0}
(-)^r[C_r]~.$$
{\rm (ii)} Isomorphism torsion defines a forgetful map
$$\begin{array}{l}
i_*~:~K^{iso}_1(\SD(\bA)) \to K^{iso}_1(\bA)~;~
\tau^{iso}(f) \mapsto [\newt(f)]~=~\newt(f)
\end{array}$$
which is a surjection split by the injection
$$K^{iso}_1(\bA) \to K^{iso}_1(\SD(\bA))~;~\tauiso(f:A \to B) \mapsto
\newt(f:(A,0) \to (B,0))~.$$
{\rm (iii)} The diagram
$$\xymatrix{
K_0(\SD(\bA)) \otimes K_0(\SD(\bA))
\ar[d]_-{\di{\chi \otimes \chi}}^-{\di{\cong}}
\ar[r]^-{\di{\epsilon}}
& K^{iso}_1(\SD(\bA)) \ar[d]^-{\di{i_*}} \\
K_0(\bA) \otimes K_0(\bA) \ar[r]^-{\di{\epsilon}} &
K^{iso}_1(\bA)}$$
commutes,
that is the sign of objects $(C,\eta_C)$, $(D,\eta_D)$ in $\SD(\bA)$
has image
$$i_*\epsilon((C,\eta_C),(D,\eta_D))~=~
\epsilon(\chi(C),\chi(D))\in K^{iso}_1(\bA)~.$$
\end{prop}
\begin{proof} (i) A short exact sequence $0 \to C \to D \to E \to 0$
of finite chain complexes in $\bA$
determines a relation
$$[C,\eta_C] - [D,\eta_D]+[E,\eta_E]~=~0 \in K_0(\SD(\bA))$$
for any signs $\eta_C,\eta_D,\eta_E$.\\
(ii) By construction.\\
(iii) The sign
$$\begin{array}{l}
\epsilon((C,\eta_C),(D,\eta_D))\\[.6ex]
=~\newt(\begin{pmatrix}
0 & 1 \\ 1 & 0 \end{pmatrix}:
(C,\eta_C)\oplus (D,\eta_D) \to (D,\eta_D) \oplus (C,\eta_C))
\in K^{iso}_1(\SD(\bA))
\end{array}$$
has image
$$\begin{array}{l}
i_*\epsilon((C,\eta_C),(D,\eta_D))\\[2ex]
\hskip25pt =~\newt(\begin{pmatrix}
0 & 1 \\ 1 & 0 \end{pmatrix}:
(C,\eta_C)\oplus (D,\eta_D) \to (D,\eta_D) \oplus (C,\eta_C))\\[2ex]
\hskip25pt =~\sum\limits^{\infty}_{r=0}(-)^r\newt
(\begin{pmatrix} 0 & 1 \\ 1 & 0 \end{pmatrix}:C_r \oplus D_r \to
D_r \oplus C_r)+\eta_{D\oplus C}-\eta_{C\oplus D}\\[2ex]
\hskip25pt =~\sum\limits^{\infty}_{r=0}(-)^r\epsilon(C_r,D_r)-
\epsilon(\chi(D),\chi(C))
+\sum\limits^{\infty}_{r=0}(-)^r\epsilon(D_r,C_r)\\[2ex]
\hskip25pt =~\epsilon(\chi(D),\chi(C))~=~\epsilon(\chi(C),\chi(D)) \in
K^{iso}_1(\bA)~.
\end{array}$$
\end{proof}

\section{Duality properties of absolute torsion.}
\label{section:dual}

In this section we extend the notion of absolute torsion to encompass
dual objects and dual maps.  We now work over an additive category
with involution
(defined below)  and introduce the notion of a dual signed
complex $C^{n-*}$ (also defined below).  We prove the following result:

\begin{prop}
\begin{enumerate}
\label{prop:dual}
\item Let $C$ be a contractible signed complex.  Then
\[\taunew(C^{n-*}) =(-)^{n+1}\taunew(C)^*\in \Kiso(\bA)\]
\item Let $f:C \to D$ be a chain equivalence of signed chain
complexes.  Then
\[\taunew(f^{n-*}:C^{n-*}\to D^{n-*}) = (-)^n\taunew(f)^*\in \Kiso(\bA)\]
\item Let $f:C^{n-*} \to D$ be a chain equivalence of signed chain
complexes.  Then the chain equivalence $Tf:D^{n-*}\to C$ (defined
below) satisfies:
\[\taunew(Tf:D^{n-*} \to C) = (-)^n\taunew(f)^* +
\frac{n}{2}(n+1)\epsilon(\chi(C),\chi(C)) \in \Kiso(\bA)\]
\end{enumerate}
\end{prop}

The rest of this section will be concerned with defining these
concepts and proving proposition \ref{prop:dual}.

Following \cite{andrew:additive} we define an \emph{involution} on
an additive category $\bA$ to be a contravariant functor
\[^*:\bA \to \bA;\;M\to M^*,\;\;(f:M\to N) \to(f^*:N^*\to M^*)\]
together with a natural equivalence
\[e:id_\bA\to **\colon \bA \to\bA;\;M\to (e(M):M \to M^{**})\]
such that for any object $M$ of $\bA$
\[e(M^*) = (e(M)^{-1})^*:M^* \to M^{***}\]

An involution on $\bA$ induces an involution on $\Kiso(\bA)$ in the
obvious way.

Throughout the rest of this chapter \(\bA\) is an additive category
with involution.

\begin{sign}
Given a \(n\)-dimensional chain complex
\[  C\;:\;C_{n} \xrightarrow{d} C_{n-1} \xrightarrow{d} C_{n-2}
\xrightarrow{d} \ldots\xrightarrow{d}C_0\]
we use the following sign convention for the dual complex \(C^{n-*}\).
\[d_{C^{n-*}} = (-)^rd_C^* \;:\;C^{n-r} \rightarrow C^{n-r+1}\]
\end{sign}
We define the sign term
\[\alpha_n(C) = \sum_{r \equiv
n+2,n+3\;(\mathrm{mod}\;4)}\epsilon(C^r,C^r)
\in \Kiso(\bA)\]

Given a signed chain complex \((C,\eta_C)\) we define the \emph{dual
signed chain complex with} \((C^{n-*},\eta_{C^{n-*}})\) by
\[\eta_{C^{n-*}}:= (-)^{n+1}\eta_C^* + (-)^{n+1}\beta(C,C)^* +
\alpha_{n}(C) \in
\Kiso(\bA)\]

\begin{lemma}
\label{lemma:dualprops}
Let \(A\) and $B$ be elements of $\bA$ and $C$ and $D$ chain complexes
over $\bA$.  We have the following basic properties of the absolute
torsion in an
additive category with involution:
\begin{enumerate}
\item \(\chi(C^{n-*}) = (-)^n\chi(C)^*\)
\item \(\epsilon(A^*,B^*) = \epsilon(B,A)^*\)
\item \(\epsilon(C^{n-*},D^{n-*}) = (-)^n\epsilon(C,D)^*\)
\item \(\beta(C^{n-*},D^{n-*}) = (-)^n\beta(D,C)^*\)
\item For chain isomorphisms $f:C\to D$ we have that:
\[\taunew(f^{n-*}:D^{n-*} \to C^{n-*}) = (-)^n\taunew(f)^* \in \Kiso(\bA)\]
\item
\label{lemma:dualprops:sum}
\[
\taunew((C\oplus D)^{n-*} \to C^{n-*}\oplus D^{n-*}) =
\left\{\begin{array}{c}0 \\
\epsilon(\chi(D),\chi(C))^*\end{array}\right.
\;\mathrm{for}\;\left\{\begin{array}{c}n\;\mathrm{even}\\
n\;\mathrm{odd}\end{array}\right.\]
\item $\taunew(1\colon C^{n-*} \to (SC)^{n+1-*}) = 0$.
\item $\taunew((-1)^{r}\colon C^{n+1-*} \to S(C^{n-*})) = 0$.
\item $\taunew((-1)^{(n+1)r}\colon (C^{n-*})^{n-*} \to C) =
\frac{n}{2}(n+1)\epsilon(\chi(C),\chi(C))^*$
\end{enumerate}
\end{lemma}
\begin{proof}
Parts 1 to 5 follow straight from the definitions.  For part 6:
\begin{eqnarray*}
&&\taunew((C\oplus D)^{n-*} \to C^{n-*}\oplus D^{n-*}) \\
&&\hspace{1cm}\begin{array}{ll}= &
\eta_{C^{n-*}\oplus D^{n-*}} - \eta_{(C\oplus D)^{n-*}} \\
 = & \epsilon(\chi(C^{n-*}),(D^{n-*})_{even}) +
(-)^n\epsilon(\chi(C),D_{even})^*
\end{array}
\end{eqnarray*}
The result follows after considering the odd and even cases.

Part 7 follows straight from the definitions.  For part 8:
\begin{eqnarray*}
&&\taunew((-1)^{r}\colon C^{n+1-*} \to S(C^{n-*})) \\
&&\hspace{1cm}\begin{array}{ll}
&\taunew((-1)^{r}\colon C^{n+1-*} \to S(C^{n-*})) \\
=&\eta_{S(C^{n-*})} - \eta_{C^{n+1-*}} + \sum_{r \equiv
n\;(\mathrm{mod}\;2)} \epsilon(C_r,C_r)^* \\
=&\alpha_{n+1}(C) + \alpha_n(C) +
\sum_{r \equiv
n\;(\mathrm{mod}\;2)} \epsilon(C_r,C_r)^* \\
=& 0
\end{array}
\end{eqnarray*}
For part 9:
\begin{eqnarray*}
&&\taunew((-1)^{(n+1)r}\colon (C^{n-*})^{n-*} \to C) \\
&&\hspace{1cm}\begin{array}{ll}
=& \eta_C - \eta_{(C^{n-*})^{n-*}} + \tau((-)^{(n+1)r}:C_r\to C_r) \\
=&\alpha_n(C^{n-*}) + (-)^n\alpha_n(C)^* +
(n+1)\sum_{r\;\mathrm{odd}}\epsilon(C_r,C_r) \\
=& \sum_{r \equiv n+2,n+3 (\mathrm{mod}\;4)}(\epsilon(C_r,C_r) +
 \epsilon(C_{n-r},C_{n-r})) \\
&+(n+1)\sum_{r\;\mathrm{odd}}\epsilon(C_r,C_r)
 \\
= & \left\{\begin{array}{c}0 \\ \epsilon(\chi(C),\chi(C)) \\
\epsilon(\chi(C),\chi(C)) \\ 0
\end{array}\right. \mathrm{for}\;n \equiv\left\{\begin{array}{c} 0
\\ 1 \\ 2 \\ 3\end{array}\right. \\
= & \frac{n}{2}(n+1)\epsilon(\chi(C),\chi(C))
\end{array}
\end{eqnarray*}
\end{proof}

\begin{lemma}
\label{lemma:dualcontract}
The torsion of a contractible signed chain
complex \(C\) in $\bA$ satisfies:
\[\taunew(C^{n-*}) = (-)^{n+1}\taunew(C)^* \in \Kiso(\bA)\]
\end{lemma}
\begin{proof}
We denote by \(\bar{C}^{n-*}\) the chain complex with
\((\bar{C}^{n-*})_r=(C_{n-r})^*\) and
\[
d_{\bar{C}^{n-*}}=d_C^*\colon\bar{C}^{n-r}\rightarrow \bar{C}^{n-r+1}
\]
We have an isomorphism \(f\colon C^{n-*}_r \rightarrow
\bar{C}^{n-*}_r\) given by \(f=-1\) if \(r \equiv
n+2,n+3\;(\mathrm{mod}\;4)\) and \(f=1\) otherwise.  By considering the
torsion of this isomorphism we have:
\begin{equation}
\label{eqn:cbar}
\tau(C^{n-*}) = \tau(\bar{C}^{n-*}) + \alpha_n(C)
\end{equation}
Let \(n_{even}\) be the greatest even integer \(\leq n\), similarly
\(n_{odd}\).
For any chain contraction \(\Gamma\) for \(C\) we have the following
commutative diagram:
\small
\[
\xymatrix{
C^{n_{even}}\oplus\ldots\oplus C^0
\ar[rrr]^*+{\left(\begin{array}{cccc} d^*&0&0&\ldots\\\Gamma^*&d^*&0&\ldots\\
0&\Gamma^*&d^*&\ldots\\\vdots&\vdots&\vdots&\end{array}\right)}
\ar[ddd]_*+{\left(\begin{array}{cccc}& & & 1\\& & 1 & \\& \ddots \\
1\end{array}\right)}
&
&
&
C^{n_{odd}}\oplus\ldots\oplus C^1
\ar[ddd]^*+{\left(\begin{array}{cccc}& & & 1\\& & 1 &
\\& \ddots \\
1\end{array}\right)} \\ \\ \\
C^0\oplus\ldots\oplus C^{n_{even}}
\ar[rrr]_*+{\left(\begin{array}{cccc}
d^*&\Gamma^*&0&\ldots\\0&d^*&\Gamma^*&\ldots\\
0&0 &d^*&\ldots\\\vdots&\vdots&\vdots&\end{array}\right)}
&&&
C^1\oplus\ldots\oplus C^{n_{odd}}
}
\]
\normalsize
The torsion of the lower map in this diagram is \(\tau(C)^*\); the
torsion of the uppermost map is \((-)^{n+1}\tau(\bar{C}^{n-*})\).  So, by first
considering the torsions of the maps in the above diagram we have:
\begin{eqnarray*}
\tau(\bar{C}^{n-*}) & = & (-)^{n+1}\tau(C)^* + (-)^{n+1}(\sum_{i>j;\; i,j\;
even}\epsilon(C^i,C^j)-\sum_{i>j;\; i,j\;
odd}\epsilon(C^i,C^j)) \\
& = &  (-)^{n+1}\tau(C)^* + (-)^n\beta(C,C)^*
\end{eqnarray*}
Hence by equation \ref{eqn:cbar}
\[\tau(C^{n-*}) = (-)^{n+1}\tau(C)^* +
(-)^n\beta(C,C)^*+\alpha_n(C)\]
Using the definition of the dual signed chain complex we have:
\[\taunew(C^{n-*}) = (-)^{n+1}\taunew(C)^*\]
\end{proof}

\begin{lemma}
\label{lemma:dualmap}
Let \(C,D\) be \(n\)-dimensional signed chain complexes
in \(\bA\) and \(f:C\rightarrow D\) a chain equivalence.  Then
\[\taunew(f^{n-*}:D^{n-*}\rightarrow C^{n-*}) =
(-)^{n}\taunew(f)^*\in \Kiso(\bA)\]
\end{lemma}
\begin{proof}
We have an isomorphism of chain complexes \(\theta:\CC(f^{n-*})\rightarrow
\CC(f)^{n+1-*}\) given by:
\[\CC(f^{n-*})_r=C^{n-r}\oplus D^{n-r+1}
\xrightarrow{\left(\begin{array}{cc} 0 & (-)^{n-r}\\ 1 &
0\end{array}\right)}\CC(f)^{n+1-*}_r= D^{n-r+1}\oplus C^{n-r}\]
The torsion of the map $\theta$ is given by:
\begin{eqnarray*}
\taunew(\theta) & = & \taunew((-)^{n-r}:S(D^{n-*}\to D^{n+1-*}) \\
&&+\newt(C^{n-*}\to (SC)^{n+1-*}) \\
&&+\newt((SC)^{n+1-*}\oplus D^{n+1-*} \to (SC \oplus D)^{n+1-*}) \\
&&+(-)^{n+1}\newt(D\oplus SC \to SC\oplus D)^* \\
&=& n\epsilon(\chi(D),\chi(D))^*+(n+1)\epsilon(\chi(D),\chi(D))^*\\
&&+\epsilon(\chi(D),\chi(D))^*\\
&=& 0
\end{eqnarray*}
and the result follows since \(\taunew(\CC(f)^{n+1-*}) = (-)^n\taunew(f)^*\)
\end{proof}

We define the duality isomorphism \(T\) as:
\[T\;:\;\mathrm{Hom}_\bA(C^p,D_q)\rightarrow
\mathrm{Hom}_\bA(D^q,C_p)\;;\;\phi \rightarrow (-)^{pq}\phi^*\]

\begin{lemma}
\label{lemma:Tdual}
Let \(C,D\) be \(n\)-dimensional signed chain complexes in $\bA$ and
\(f:C^{n-*} \rightarrow D\) a chain equivalence.  Then
\[\taunew(Tf:D^{n-*} \to C) = (-)^n\taunew(f)^* +
\frac{n}{2}(n+1)\epsilon(\chi(C),\chi(C)) \in \Kiso(\bA)\]
\end{lemma}
\begin{proof}
Using lemmas \ref{lemma:dualprops} and \ref{lemma:dualmap} we have:
\begin{eqnarray*}
\taunew(Tf:D^{n-*}\to C) &=& \taunew(f^*:D^{n-*} \to(C^{n-*})^{n-*})
\\
&&+
\taunew((-1)^{(n+1)r}(C^{n-*})^{n-*} \to C) \\
& = & (-)^n\taunew(f)^*+
\frac{n}{2}(n+1)\epsilon(\chi(C),\chi(C))
\end{eqnarray*}
as required.
\end{proof}

Together the above three lemma prove proposition \ref{prop:dual}.

\section{Torsion of Poincar\'{e} complexes}
\label{section:poincare}

We now move on to consider symmetric Poincar\'e complexes.  These are
algebraic objects which encapsulate the properties of Poincar\'e
duality spaces (see Ranicki \cite{origandrew},
\cite{origandrew2} for a more complete discussion).
We will restrict ourselves to considering symmetric Poincar\'e
complexes over a ring $R$, that is we work over $\bA=\bA(R)$ and we
will consider the torsion invariants to lie in the more familiar
$K_1(R)$.
 We will define
the notion of the absolute torsion
of
a symmetric Poincar\'e complex to be, essentially, the torsion of the
Poincar\'e duality chain equivalence.  In the case of compact oriented
manifolds $M^n$ with $CW$-structure this is the torsion of the map:
\[[M^n] \cap -: C(M)^{n-*} \to C(M)\]
In this section we develop
the theory from this algebraic viewpoint; it will be applied to
geometric objects in a later section.

We recall from \cite{origandrew} and \cite{aar_topological} the following
definition:
\begin{defin}
\begin{enumerate}
\item An \emph{$n$-dimensional symmetric complex}
\((C,\phi_0)\) is  a chain complex \(C\) in $\bA(R)$, together
with a collection of  morphisms
\[\phi = \{\phi_s\colon C^{n-r+s} \rightarrow C_r\;|\;s\geq 0\}\]
such that
\begin{eqnarray*}
& & d_c\phi_s +
(-)^r\phi_sd^*_C+(-)^{n+s+1}(\phi_{s-1}+(-)^sT\phi_{s-1}) = 0 \\
& & \colon
C^{n-r+s-1}\rightarrow C_r\;\;
 (s\geq 0, \phi_{-1}=0)
\end{eqnarray*}
Hence \(\phi_0\) is a chain map and \(\phi_1\) is a chain homotopy
\(\phi_1: \phi_0 \simeq T\phi_0\).
\item The complex is said to be \emph{Poincar\'e} if \(\phi_0\) is a
chain equivalence.
\item The complex is said to be \emph{round} or \emph{even} if \(C\) is round
or even respectively.
\item A \emph{morphism} between \(n\)-dimensional symmetric complexes
\((C,\phi)\) and \((C',\phi')\) consists of a chain map \(f:C
\rightarrow C'\) and morphisms \(\sigma_s: {C'}^{n+1+s-r} \rightarrow
C'_r\) \(s \geq 0\) such that
\[\phi'_s - f\phi_sf^* =
d\sigma_s+(-)^r\sigma_sd^*+(-)^{n+s}(\sigma_{s-1}+(-)^sT\sigma_{s-1}):C^{n-r+s}
\rightarrow C_r\]
(in particular \(\phi_0' \simeq f\phi_0f^*\)).  Such a morphism is
said to be a \emph{homotopy equivalence} if \(f\) is a chain equivalence.
\item A symmetric complex \((C,\phi)\) is said to be \emph{connected}
if \(H_0(\phi_0:C^{n-*} \rightarrow C)=0\)
\item The \emph{boundary} \((\partial C,\partial\phi)\) of a connected
\(n\)-dimensional symmetric complex \((C,\phi)\) is the
\((n-1)\)-dimensional symmetric Poincar\'e complex defined by
\begin{eqnarray*}
&d_{\partial C} = &\left(\begin{array}{cc}d_C & (-)^r\phi_0 \\ 0 &
(-)^rd^*_C \end{array}\right):\partial C_r = C_{r+1} \oplus C^{n-r}
\rightarrow \\
&&\partial C_{r-1} = C_r \oplus C^{n+1-r} \\
&\partial\phi_0 = &\left(\begin{array}{cc} (-)^{n-r-1}T\phi_1 &
(-)^{rn} \\ 1 & 0\end{array}\right):\partial C^{n-r-1} = C^{n-r}
\oplus C_{r+1} \rightarrow \\
&&\partial C_r = C_{r+1} \oplus C^{n-r} \\
&\partial\phi_s = &\left(\begin{array}{cc} (-)^{n-r-1}T\phi_{s+1} &
(-)^{rn} \\ 0 & 0\end{array}\right):\partial C^{n-r+s-1} = C^{n+s-r}
\oplus C_{r-s+1} \\
&& \rightarrow\partial C_r = C_{r+1} \oplus C^{n-r}
\end{eqnarray*}
\item A signed symmetric (Poincar\'e) complex is a
symmetric (Poincar\'e) complex \((C,\phi_0)\) where in addition
\(C\) is a signed chain complex.
\end{enumerate}
\end{defin}

\begin{example}  An $n$-dimension manifold $M$ with universal covering
$\widetilde{M}$ determines a symmetric Poincar\'e complex
$(C(\widetilde{M}),\phi)$ in $\bA(\Z[\pi_1 M])$ with
\[\phi_0 = [M] \cap -:C(\widetilde{M})^{n-*} \to C(\widetilde{M})\]
\end{example}

\begin{lemma}
\label{prop:boundary}
The boundary $(\partial C,\partial\phi)$ of any signed $n$-dimensional
symmetric
complex \((C,\phi)\) satisfies
\[\taunew(\partial\phi_0:(\partial C)^{n-1-*} \to \partial C) =
\frac{n}{2}(n+1)\epsilon(\chi(C),\chi(C))
\in K_1(R)\]
\end{lemma}
\begin{proof}
The map
\[\partial\phi_0 = \left(\begin{array}{cc} (-)^{n-r-1}T\phi_1 &
(-)^{rn} \\ 1 & 0\end{array}\right):\partial C^{n-r-1}
\rightarrow \partial C_r\]
is an isomorphism.

We have that
\begin{eqnarray*}
\taunew(\partial\phi_0)&=&\taunew\left(\left(\begin{array}{cc}
0&(-)^rn\\1&0\end{array}\right):(\Omega C\oplus C^{n-*})^{n-1-*} \to
\Omega C \oplus C^{n-*}\right)\\
&=&\taunew((\Omega C\oplus C^{n-*})^{n-1-*} \to (\Omega
C)^{n-1-*}\oplus (C^{n-*})^{n-1-*}) \\
&&+ \taunew((C^{n-*})^{n-1-*} \to
(\Omega C^{n-1-*})^{n-1-*}) \\
&&+ \taunew((-)^{nr}:(\Omega
C^{n-1-*})^{n-1-*} \to \Omega C)\\
&& + \taunew((\Omega C)^{n-1-*} \to
C^{n-*})\\
&=&\taunew(C^{n-*}\oplus \Omega C\to \Omega C \oplus C^{n-*})\\
&=& \frac{n}{2}(n+1)\epsilon(\chi(C),\chi(C))
\end{eqnarray*}
using the results of lemma \ref{lemma:dualprops}.
\end{proof}

We can now define a new absolute torsion invariant of Poincar\'{e}
complexes which is additive and a cobordism invariant.

\begin{defin}
\label{defin:newtorsion}
We define the \emph{absolute torsion} of a signed Poincar\'{e}
complex \((C,\phi)\) as
\[\tau^{NEW}(C,\phi) = \taunew(\phi_0) \in K_1(R)\]
\end{defin}
\begin{prop}
\label{prop:newtorsion}
Let \((C,\phi)\) and \((C',\phi')\) be signed $n$-dimensional
Poincar\'{e} complexes. Then:
\begin{enumerate}
\item
\label{prop:newtorsion:additive}
\emph{Additivity:}
\label{prop:newtorsion:additivity}
\[\tau^{NEW}(C\oplus C',\phi \oplus \phi') =
\tau^{NEW}(C,\phi) + \tau^{NEW}(C'\phi_0') \in K_1(R)\]
\item \emph{Duality:}
\label{prop:newtorsion:duality}
\[\tau^{NEW}(C,\phi) =
(-)^n\tau^{NEW}(C,\phi)^*+
\frac{n}{2}(n+1)\epsilon(\chi(C),\chi(C))\in K_1(R)\]
(n.b. the above sign term disappears in the case where
$\bA=\bA(R)$ and where
anti-symmetric
form over the ring \(R\) necessarily have even rank; this is the case
for \(R=\Z\) or \(R=\Q\) but not \(R=\C\)).
\item \emph{Homotopy invariance:}
\label{prop:newtorsion:homotopy}
Suppose \((f,\sigma_s)\) is a
homotopy equivalence from \((C,\phi)\) to \((C',\phi')\). Then
\[\taunew(C',\phi') = \taunew(C,\phi) + \tau(f) + (-)^n\tau(f)^*\in
K_1(R)\]
\item
\label{prop:newtorsion:cobordism}
\emph{Cobordism Invariance:} Suppose that \((C,\phi)\) is
homotopy equivalent to the boundary of some \((n+1)\)-dimensional
symmetric complex with torsion
\((D,\phi^D)\).  Then
\begin{eqnarray*}
\tau^{NEW}(C,\phi)&=&(-)^{n+1}\taunew(C \rightarrow \partial D)^*-\taunew(C
\rightarrow \partial D) \\ &&+
\frac{1}{2}(n+1)(n+2)\epsilon(\chi(D),\chi(D)) \in K_1(R)
\end{eqnarray*}
\item \emph{Orientation change:}
\[\tau^{NEW}(C,-\phi) = \tau^{NEW}(C,\phi) + \epsilon(\chi(C),\chi(C))
\in K_1(R)\]
\item
\label{prop:newtorsion:torsion}
The absolute torsion of a
signed Poincar\'e complex is independent of
the choice of sign \(\eta_C\).
\end{enumerate}
\end{prop}
\begin{proof}
\begin{enumerate}
\item
A symmetric Poincar\'e complex of odd dimension satisfies $\chi(C)=0$,
hence the map $(C\oplus C')^{n-*}\to C^{n-*} \oplus {C'}^{n-*}$ has
trivial absolute torsion.  Additivity now follows from the additivity
of chain equivalences.
\item We know that \(\phi_0\) is homotopic to \(T\phi_0\); duality now
follows by applying lemma \ref{lemma:Tdual}.
\item We have that \(\phi_0' \simeq f\phi_0f^*\) and hence
\[
\taunew(\phi'_0) = \taunew(f) + \taunew(\phi_0) + (-)^n\taunew(f)^*
\]
\item This follows from proposition \ref{prop:boundary} and
homotopy invariance.
\item We have that \(\taunew(-\phi_0) = \taunew(\phi_0) + \taunew(-1:C
\rightarrow
C) = \taunew(\phi_0) + \chi(C)\tau(-1)\).
\item A change in \(\eta_C\) leads to a corresponding change in
\(\eta_{c^{n-*}}\) so \(\taunew(\phi_0)\) is unchanged.
\end{enumerate}
\end{proof}

\section{The signed Poincar\'e derived category with involution.}
\label{section:spd}

In this section we will add an involution to a particular subcategory
of the signed derived category.  Let $\SPD_n(\bA(R))$ denote the category
whose object are signed $n$-dimensional chain complexes $C$ in
$\bA(R)$ which are isomorphic to their dual complexes $C^{n-*}$ and
$\chi(C)=0$ if$n$ is odd.  Then
we have an involution
\[*: C \mapsto C^{n-*}\]
\[*:(f:C\to D) \mapsto (f^{n-*}:D^{n-*}\to C^{n-*})\]
with the natural equivalence $e(C)$ given by
\[e(C)=(-)^{(n+1)r}:C\to (C^{n-*})^{n-*}\]
We call this category the \emph{signed Poincar\'e derived category
  with $n$-involution}.
In order to show that this is a covariant functor of additive
categories we must show that $*(A\oplus B)=*A\oplus *B$.  However, the
  condition that $\chi(C)$ is odd if $n$ is odd implies that the
  torsion of the rearrangement map $(C\oplus
D)^{n-*}\to C^{n-*}\oplus D^{n-*}$ is trivial (see lemma
  \ref{lemma:dualprops} part \ref{lemma:dualprops:sum}) and the functor $*$
  is additive.
As in the case of the signed derived category we have a map
\[i_*:\Kiso(\SPD_n(\bA(R))) \to \Kiso(\bA(R))\;;\;\tau^{iso}(f)\mapsto
\taunew(f)\]
The behaviour of $i_*$ under the involution on $\SPD_n(\bA(R))$ is given by
\[i_*(f^*) = (-)^ni_*(f)^*\]

\section{A product formula}
\label{section:product}

In this section we will quote a formula for the absolute torsion of a
product of symmetric Poincar\'e complexes and prove it in a special case.
We recall from \cite{origandrew} the definition of the tensor product
of symmetric Poincar\'e complexes:

\begin{defin}
\begin{enumerate}
\item
The tensor product $C\otimes D$ of a chain complex $C$ in $\bA(R)$ and
a
chain complex $D$ in $\bA(R')$ is the chain complex in $\bA(R\otimes R')$
\[d_{C \otimes D} : (C\otimes D)_r =\sum_{s=-\infinity}^\infinity
C_s\otimes D_{r-s}\to
(C\otimes D)_{r-1};\]
\[x \otimes y \mapsto x \otimes d_D(y)+(-)^{r-s}d_C(x) \otimes y\]
\item The tensor product $(C\otimes D,\phi \otimes \theta)$ of an
  $n$-dimensional symmetric Poincar\'e complex $(C,\phi)$ with an
  $m$-dimensional symmetric Poincar\'e complex $(D,\theta)$ is an
  $(n+m)$-dimensional Poincar\'e complex defined by:
\[(\phi\otimes\theta)_s = \sum_{r=0}^s(-1)^{(n+r)s}\phi_r\otimes
  T\theta_{s-r} : (C\otimes D)^{n+m-*} \to C\otimes D\]
\end{enumerate}
\end{defin}

\begin{prop}
Let $(C,\phi)$ and $(C',\phi')$ be symmetric Poincar\'e complexes over
a rings with involution $R$ and $R'$ respectively.  Then:

\[\taunew(C\otimes C',\phi \otimes \phi') = \chi(C)\taunew(C',\phi') +
\chi(C')\taunew(C,\phi)\]
\end{prop}

The proof of this result requires the theory of signed complexes to be
extended to tensor products; this theory is developed in \cite{hkr}
using the theory of signed filtered complexes.

We will use the following ``ad hoc'' methods to
prove the product formula under the following condition:
\begin{assumption}
\label{assumption:rings}
The rings $R$ and $R'$ are such that $4k+2$-dimensional symmetric
forms necessarily have even
rank (e.g. group rings).
\end{assumption}
The ring $R=\C$ is an example which does not satisfy this assumption.
We say a module $M$ in $\bA(R)$ is \emph{even} if $\rR(M)$ is even;
similarly we say a chain complex $C$ in $\bA(R)$ is
\emph{even} if $C_r$ is even for all $r$.
\begin{defin}
Let $(C,\eta_C)$ and $(D,\eta_D)$ be even signed complexes.  We define the
\emph{signed complex tensor product} by
\[(C\otimes D,\eta_{C\otimes D}) = (C\otimes D,
\chi(C)\eta_{D}+\chi(D)\eta_{C})\]
\end{defin}

\begin{lemma}
\label{lemma:products}
\begin{enumerate}
\item Let $C$, $C'$ be signed, even complexes in $\bA(R)$ and $D$ a
  signed even complex in $\bA(R')$.  Then
\[\taunew((C\oplus C')\otimes D \to (C\otimes D)\oplus(C'\otimes
  D))=0\]
\[\taunew(D\otimes(C\oplus C')
\to (D\otimes C)\oplus(D\otimes C'))=0\]
\item Let $C$ be a signed contractible even complex in $\bA(R)$
  and $D$ a signed even complex in $\bA(R')$.  Then:
\[\taunew(C\otimes D) = \taunew(D \otimes C) = \chi(D)\taunew(C)\]
\item Let $f:C \to C'$ be a chain equivalence of even complexes
  in $\bA(R)$ and $D$ a signed even complex in $\bA(R')$.  Then:
\begin{eqnarray*}
\taunew(f \otimes 1:C\otimes D \to C'\otimes D)&=&\taunew(1\otimes
  f:D \otimes C \to D \otimes C') \\
&=& \chi(D)\taunew(f)
\end{eqnarray*}
\end{enumerate}
\end{lemma}
\begin{proof}
For even modules $M,N$ the map $\epsilon(M,N)=0 \in K_1(R)$ and
$\tau(-1:M \to M)=0$, so the torsion of rearrangement maps is always
zero and $\eta_{C\oplus D} =
\eta_D + \eta_C$ for even chain complexes $C$, $D$.  Part 1 follows straight
from these facts.  For part 2, let
$\Gamma:C_* \to C_{*+1}$ be a chain contraction of $C$,  then
$(d_C+\Gamma)\otimes 1$ is a chain contraction of $C\otimes D$. We have
a commutative diagram:
\[\xymatrix{(C\otimes D)_{odd} \ar[r]^{(d_C+\Gamma)\otimes 1} \ar[d]
  & (C\otimes
  D)_{even}  \ar[d] \\
(C_{even}\otimes D_{odd})\oplus(C_{odd}\oplus D_{even})
\ar[r]
& (C_{odd}\otimes D_{odd})\oplus(C_{even}\oplus D_{even})
}\]
with the bottom map given by
$(d_C+\Gamma)\otimes 1 \oplus (d_C+\Gamma)\otimes 1$.
The torsion of the top map is $\tau(C\otimes D)$, the torsion of the
bottom map is $\chi(D)\tau(C)$.  Part 2 follows from the fact that
the torsions of the left and right maps are zero, since they are
rearrangements.  For part 3 we have that
\[\taunew(\CC(f\otimes 1) \to \CC(f)\otimes D)=0\]
since it is a rearrangement map and the result now follows straight
from the definitions.
\end{proof}

\begin{prop}
Let \((C,\phi^C)\) and \((D,\phi^D)\) be Poincar\'e complexes over
ring $R$ and $R'$ respectively, where $R$ and $R'$ satisfy assumption
\ref{assumption:rings}.  Then
\begin{eqnarray*}
\taunew(C \otimes D,\phi^C \otimes \phi^D)&=&\chi(D)\taunew(C,\phi^C) +
\chi(C)\taunew(D,\phi^D) \\
&&\in K_1(R \otimes R')
\end{eqnarray*}
\end{prop}
\begin{proof}
Given any f.g. based chain complex \(C\) we may form a direct sum
with a contractible chain complex to form a new chain complex which is
of even rank in every dimension except one.  Hence given a Poincar\'e
pair \((C,\phi)\) we may form a new Poincar\'e pair \((C',\phi^{C'})\) which
is even in every dimension by
\begin{enumerate}
\item Forming the direct sum with a contractible complex (letting
\(\phi^{C'}\) vanish on this contractible complex) such that \(C'\) is
even in every dimension except possibly the middle.
\item If the complex is odd (and hence of dimension 4k), forming the
direct sum with the
Poincar\'e complex which is \(R\) in dimension \(2k\), vanishes
otherwise and has \(\phi_0 = 1:R \to R\).
\end{enumerate}
We may form a similar complex \((D',\phi^{D'})\) from \((D,\phi^D)\).
Using lemma \ref{lemma:products} we see that
\[\taunew(\phi^{C'}_0 \otimes
\phi^{D'}_0)=\chi(C')\taunew(\phi_0^{D'}) +
\chi(D')\taunew(\phi_0^{C'}) \in K_1(R\otimes R')\]
Hence
\[\taunew(C' \otimes D',\phi^{C'} \otimes
\phi^{D'})=\chi(C')\taunew(D',\phi_0^{D'}) + \chi(D')\taunew(C',\phi_0^{C'})
\]
Let \(R^C\) denote the chain complex which is 0 if \(C\) is even and
\(R\) in dimension \(2k\) otherwise; similarly \({R'}^D\).
By direct computation (for the sake of clarity we now suppress mention
of the morphisms \(\phi\))
\[\taunew(R^C \otimes D) = \chi(C)(\taunew(D) +
(k+1)\chi(D))\tau(-1))\]
\[\taunew(C \otimes {R'}^D) = \chi(D)(\taunew(C) + (l+1)\chi(C)\tau(-1))\]
\[\taunew(R^C \otimes {R'}^D) = (k+l)\chi(C)\chi(D)\tau(-1)\]
The Poincar\'e complex \((C' \otimes D' )\) is homotopy equivalent to
\((C\oplus R^C)\otimes(D\oplus{R'}^D)\); using the invariance of
absolute torsion under homotopy equivalence, its additivity
properties and the above three formulae we have:
\[\taunew(C \otimes D) = \chi(C)\taunew(D) + \chi(D)\taunew(C) \in
K_1(R\otimes R')\]
as required.
\end{proof}

\section{Round L-theory}
\label{section:rltheory}

We refer the reader to \cite{rltheory} for the definition of the round
symmetric $L$-groups \(L^n_r(A)\).  The absolute torsion defined in this
paper as \(\tau(C,\phi_0) = \tau(\phi_0)\) (here \(\tau(\phi_0)\)
refers to the absolute torsion defined in \cite{alg}) is not a cobordism
invariant.  We can define such an invariant using the absolute torsion
of a Poincar\'e complex.  If this invariant is substituted for
$\tau(C,\phi_0)$ as defined in \cite{rltheory} then the results become
correct.
\begin{lemma}
Let \((C,\phi)\) be a round Poincar\'e complex.  The reduced element
\[\taunew(C,\phi) \in \widehat{H}^n(\Z_2;K_1(R))\]
is independent of the choice of sign $\eta_C$; moreover
we have a well defined homomorphism:
\[L^n_r(R) \rightarrow \widehat{H}^n(\Z_2;K_1(R))\]
given by \((C,\phi) \mapsto \taunew(C,\phi)\).
\end{lemma}
\begin{proof}
The element \(\taunew(C,\phi) \in \widehat{H}^n(\Z_2;K_1(R))\) is
independent of the choice of sign by proposition
\ref{prop:newtorsion} part \ref{prop:newtorsion:torsion}.
The absolute torsion is additive by proposition
\ref{prop:newtorsion} part \ref{prop:newtorsion:additive}.  The
absolute torsion of the boundary of a round symmetric complex is
trivial in the reduced group \(\widehat{H}^n(\Z_1;K_1(R))\) by proposition
\ref{prop:newtorsion} part \ref{prop:newtorsion:cobordism}.  Hence the
torsion of a round null-cobordant complex is trivial and the map
\[L^n_r(R) \rightarrow \widehat{H}^n(\Z_2;K_1(R))\]
given by \((C,\phi) \mapsto \taunew(C,\phi)\) is well defined.
\end{proof}


\section{Applications to manifolds}
\label{section:manifolds}
\subsection{The absolute torsion of oriented manifolds}

To any \(n\)-dimensional oriented manifold \(M\) we may associate a
Poincar\'e complex \((C,\phi)\) over the ring \(R=\Z[\pi_1M]\) , well
defined up to homotopy
equivalence (see \cite{origandrew2}).
By property \ref{prop:newtorsion:duality} of proposition
\ref{prop:newtorsion} the absolute torsion of such a Poincar\'e
complex satisfies
\[\taunew(C,\phi) = (-)^n\taunew(C,\phi)^*\]
since $\chi(M)\equiv 0$ (mod 2) unless $n\equiv 0$ (mod 4).  Hence the torsion
\(\taunew(C,\phi)\) may be considered to lie in the group
\(\widehat{H}^n(\Z_2;K_1(\Z[\pi_1M^n]))\).  By property
\ref{prop:newtorsion:homotopy} of proposition \ref{prop:newtorsion}
if $(C,\phi)$ is homotopy equivalent to $(C',\phi')$ then
\[\taunew(C,\phi) =\taunew(C',\phi')\in \widehat{H}^n(\Z_2;K_1(\Z[\pi_1M^n]))\]
hence
\[\taunew(M^n) := \taunew(C,\phi) \in \widehat{H}^n(\Z_2; \Z[\pi_1M^n])\]
is well defined.

\subsection{Examples of the absolute torsion of manifolds.}
\subsubsection{The circle}

We may associate to the circle (\(S^1)\) the following chain complex over
$R=\Z[\pi_1(S^1)]=\Z[t,t^{-1}]$ by giving it the CW-decomposition
consisting of one 1-cell and one 0-cell:
\[
\xymatrix{
\Z[t,t^{-1}] \ar[r]^1 \ar[d]_{t^{-1}-1} &
\Z[t,t^{-1}] \ar[d]^{1-t} \\
\Z[t,t^{-1}] \ar[r]^t &
\Z[t,t^{-1}]
}
\]
In this diagram the two modules on the right are the chain complex,
the two modules on the left are the dual complex and the sideways
arrows represent \(\phi_0\).

Hence \(\taunew(S^1) = \tau(-t) \in \widehat{H}^n(\Z_2;K_1(\Z[t,t^{-1}])\).

\subsubsection{The absolute torsion of an algebraic mapping torus.}
The mapping torus of a map $f:M \to M$ is space
obtained from $M \times I$ obtained by attaching the boundaries $M
\times \{0\}$ and $M \times \{1\}$ using the map $f$.  The following
algebraic analogue is defined by Ranicki (\cite{knotty_but_nice}
definition 24.3, the reader should note the different sign convention
used here).
\begin{defin}
The \emph{algebraic mapping torus} of a morphism
$(f,\sigma):(C,\phi) \to (C,\phi)$ from an $n$-dimensional symmetric
Poincar\'e complex $(C,\phi)$ over a ring $R$ to itself is the
$(n+1)$-dimensional
symmetric complex $(T(f),\theta)$ over the ring $R[z,z^{-1}]$ defined by:
\[T(f)=\CC(f-z)\]
\[\theta_0 = \left(\begin{array}{ll}(-)^n\sigma_0 & \phi_0z \\
(-)^{n-r+1}\phi_0f^* & 0\end{array}\right):T(f)^{n+r+1} \to T(f)\]
The complex is Poincar\'e if the morphism $f$ is a chain equivalence.
\end{defin}

\begin{lemma}
\label{lemma:mapping_torus}
Let $(f,\sigma):(C,\phi) \to (C,\phi)$ be a self chain equivalence
from an $n$-dimensional symmetric
Poincar\'e complex $(C,\phi)$ over a ring $R$ to itself.  Then:
\[\taunew(T(f),\theta) = \taunew(f) + \taunew(-z:C \to C) \in
K_1(R[z,z^{-1}])\]
\end{lemma}

\begin{proof}
We have a commutative diagram with short exact rows:
\[\xymatrix@C+20pt{
0 \ar[r] & C^{n-r} \ar[r]^-{\left(\begin{array}{l}0\\1\end{array}\right)}
\ar[d]^{\displaystyle{-z\phi_0}} & \CC(f-z)^{n+1-r}
\ar[r]^-{\left(\begin{array}{ll}(-)^r&0\end{array}\right)}
\ar[d]^{\displaystyle{\theta}} & S(C^{n-r})
\ar[d]^{\displaystyle{(-)^n\phi_0f^*}} \ar[r] & 0\\
0 \ar[r] & C_r \ar[r]_-{\left(\begin{array}{l}1\\0\end{array}\right)}
& \CC(f-z)_r \ar[r]_-{\left(\begin{array}{ll}0&1\end{array}\right)} &
SC_r \ar[r] & 0
}\]
The absolute torsion of the lower short exact sequence is trivial;
for the top map we have:
\begin{eqnarray*}
\taunew\left(C^{n-r},\CC(f-z)^{n+1-r},S(C^{n-r});
\left(\begin{array}{c}0\\1\end{array}\right),
\left(\begin{array}{cc}(-)^r&0\end{array}\right)\right) = \\
\hspace{1cm}\begin{array}{l}
\taunew\left(\left(\begin{array}{cc}0&(-)^r\\1&0\end{array}\right)
:(C\oplus SC)^{n+1-*} \to C^{n-*}\oplus S(C^{n-*})\right) \\
\begin{array}{ll}
=&\taunew((C\oplus SC)^{n+1-*} \to C^{n+1-*}\oplus (SC)^{n+1-*})\\
&+\taunew((-1)^r:C^{n+1-*} \to S(C^{n-*}))\\
&+\taunew((SC)^{n+1-*}\to C^{n-*})\\
&+\taunew(S(C^{n-*})\oplus C^{n-*} \to C^{n-*} \oplus S(C^{n-*}))
\end{array} \\
=0
\end{array}
\end{eqnarray*}
Using proposition \ref{prop:torsionprops} part
\ref{prop:torsionprops:shortexact} we have:
\begin{eqnarray*}
\taunew(T(f),\theta) &=&\taunew(\theta) \\
&=&\taunew(-z\phi_0)+\taunew((-1)^nS(\phi_0f^*):S(C^{n-*}) \to SC)\\
&=&-\taunew(f^*) + \taunew(-z:C \to C)\\
&=&(-)^{n+1-*}\taunew(f)^*+\taunew(-z:C \to C)\\
&=&\taunew(f)+\taunew(-z:C \to C)
\end{eqnarray*}
as required.
\end{proof}

\subsubsection{A specific example of a mapping torus.}
We return to the example of the orientation preserving
self-homeomorphism
$f:\CP^2\to \CP^2$ given by complex conjugation in some choice of
homogeneous
coordinates (see section \ref{cp2_self_map}).  We recall that the
torsion of this map is $\taunew(f)=\tau(-1) \in K_1(\Z)$. Using lemma
\ref{lemma:mapping_torus} we compute the absolute torsion of the
mapping torus of $f$ as
\[\taunew(T(f)) = \tau(z^3) \in K_1(\Z[z,z^{-1}])\]
where $z$ is a generator of $\pi_1(T(f)) = \pi_1(S^1) = \Z$.  By
contrast we may compute the absolute torsion of the space $T(Id:\CP^2 \to
\CP^2) = S^1\times\CP^2$ as
\[\taunew(S^1\times\CP^2)=\tau(-z^3) \in K_1(\Z[z,z^{-1}])\]
hence the absolute torsion can distinguish between these two $\CP^2$
bundles over $S^1$.  A more thorough investigation into the absolute torsion
of fibre bundles of compact manifolds is made in \cite{hkr}.

\section{Identifying the sign term}
\label{section:sign}

Throughout this section we work over a group ring $R=\Z[\pi]$ for some
group $\pi$ (or, more generally, any ring with involution $R$ which
admits a map $R \to \Z$ such that the composition $\Z\to R\to\Z$ is
the identity).  We first identify the relationship between the ``sign''
term of the absolute torsion of a Poincar\'e complex and the
traditional signature and Euler characteristic and semi-characteristic
invariants.

We have a canonical decomposition of $K_1(\Z[\pi])$ as follows:
\[K_1(\Z[\pi]) = \widetilde{K_1}(\Z[\pi])\oplus \Z_2\]
with the $\Z_2$ component the ``sign'' term identified by the map
\[i_*K_1(\Z[\pi]) \to K_1(\Z)=\Z_2\]
induced by the augmentation map $i:\pi \to 1$ (more generally, a map
$R\to \Z$ gives a map $K_1(R) \to K_1(\Z)=\Z_2$ which gives a
splitting $K_1(R) = \widetilde{K_1}(R)\oplus\Z_2$).
We wish to determine the $\Z_2$ component in terms of more traditional
invariants of Poincar\'e complexes.  The augmentation map may also be
applied to
a symmetric complex $(C,\phi)$ over $\Z[\pi]$ to form a symmetric
complex over $\Z$ by forgetting the group.  Functoriality of the
absolute torsion tells us that this complex has the same sign term as
$(C,\phi)$, hence to identify the sign term it is sufficient to
consider symmetric Poincar\'e complexes over $\Z$.  We will require
the Euler semi-characteristic $\chi_{1/2}(C)$ of Kervaire \cite{kervaire}

\begin{defin}
The \emph{Euler semi-characteristic} $\chi_{1/2}(C)$ of a
$(2k-1)$-dimensional chain complex
$C$ over a field $F$ is defined by
\[\chi_{1/2}(C)=\sum_{i=0}^{k-1}(-)^i\mathrm{rank}_FH_i(C)\in \Z\]
For a $(2k-1)$-dimensional chain complex
$C$ over $\Z$ we define
\[\chi_{1/2}(C;F)=\chi_{1/2}(C\otimes_\Z F)\]
\end{defin}

\begin{prop}
\label{prop:poincaresign}
The absolute torsion of an $n$-dimensional symmetric Poincar\'e
complex over $\Z$ is determined by the signature and the Euler
characteristic and semi-characteristic as follows:
\begin{enumerate}
\item If $n=4k$ then:
\[\taunew(C,\phi)=\frac{\sigma(C)-(1+2k)\chi(C)}{2}\tau(-1)\]
with $\sigma(C)$ the signature of the complex.
\item If $n=4k+1$ then $\taunew(C,\phi) = \chi_{1/2}(C;\Q)$.
\item Otherwise $\taunew(C,\phi)=0$.
\end{enumerate}
\end{prop}

As an example we have a simple corollary:
\begin{corollary}
\label{cor:selfmaps}
The absolute torsion of an orientation preserving self-homeomorphism of
a simply-connect manifold of dimension $4k+2$ is trivial.
\end{corollary}
\begin{proof}
Let $f:M \to M$ be such a self-homeomorphism.  We may triangulate $M$
and hence construct the algebraic mapping torus $T(f)$ of the
chain equivalence $f:C(M) \to C(M)$.  By lemma
\ref{lemma:mapping_torus}
\[\taunew(T(f)) = \taunew(f) \in K_1(\Z[z,z^{-1}])\]
The augmentation map $\epsilon:\Z \to {1}$ induces a map of rings
$\epsilon_*:\Z[z,z^{-1}] \to \Z$.  Since $L$-theory and the absolute
torsion are
functorial, $\epsilon_*T(f)$ represents an element of $L^{4k+3}(\Z)$
with absolute torsion $\taunew(\epsilon_*T(f))=\taunew(f)\in K_1(\Z)$.
However
by part 3 of the above proposition $\taunew(\epsilon_*T(f))=0$.
\end{proof}

The aim of the rest of this section is to prove proposition
\ref{prop:poincaresign} .
We recall from \cite{origandrew2} the computation of the symmetric $L$-groups
$L^n_h(\Z)$ of the integers $\Z$:
\[L^n_h(\Z) =
\left\{\begin{array}{l}\Z\;(\mathrm{signature})\\\Z_2\;(\mathrm{de
  \;Rham\;invariant})\\0\\0\end{array}\right.\mathrm{for}\;n\equiv
\left\{\begin{array}{l}0\\1\\2\\3\end{array}\right.(\mathrm{mod}\;4)\]
The deRham invariant $d(C) \in Z_2$ of a $(4k+1)$-dimensional
  Poincar\'e complex was expressed in \cite{lmp} as the difference
\[d(C)=\chi_{1/2}(C;\Z_2) - \chi_{1/2}(C;\Q)\]

For dimensions $n\equiv2,3\;(\mathrm{mod}\;4)$ the absolute torsion is
a cobordism invariant (proposition \ref{prop:newtorsion} part
\ref{prop:newtorsion:cobordism}) so the above computation of the
symmetric $L$-groups tells us that the absolute torsion is trivial in
these cases, thus proving the third part of proposition
  \ref{prop:poincaresign}.

If $n=4k+1$ then the absolute torsion is not a cobordism invariant;
however it is a \textit{round} cobordism invariant, so absolute torsion
defines a map:
\[L_{rh}^{4k+1}(\Z) \to K_1(\Z)\]
Since $\chi(C)=0$ for all odd-dimensional symmetric Poincar\'e
complexes every such $(4k+1)$-dimensional complex represents an element in
$L_{rh}^{4k+1}(\Z)$.  In \cite{rltheory} (proposition 4.2) the group
$L_{rh}^{4k+1}(\Z)$ is identified as:
\[L_{rh}^{4k+1}(\Z)=\Z_2\oplus \Z_2;\;\;\;C\mapsto
(\chi_{1/2}(C;\Z_2), \chi_{1/2}(C,\Q))\]
We now construct explicit generators of this group and compute their
absolute torsions.  We define the generator $(G,\phi)$ to have chain
complex $G$
concentrated in dimensions $2k$ and $2k+1$ defined by:
\[d_{G}=0:G_{2k+1}=\Z \to G_{2k}= \Z\]
with the morphisms $\phi$ given by:
\[\phi_0=\left\{\begin{array}{l}1:G^{2k}=\Z \to G_{2k+1}=\Z\\1:G^{2k+1}=\Z
\to G_{2k}=\Z\end{array}\right.\;\;\;\phi_1=0\]
Geometrically $(G,\phi)$ is the symmetric Poincar\'e complex over $\Z$
associated to the circle.
By direct computation, $\chi_{1/2}(G,\Z_2)=1$, $\chi_{1/2}(G,\Q)=1$
 and $\taunew(G,\phi) = \tau(-1)$.
We define the generator $(H,\psi)$ to have chain complex $H$
concentrated in dimensions $2k$ and $2k+1$ defined by:
\[d_{H}=2:H_{2k+1}=\Z \to H_{2k} = \Z\]
with the morphisms $\psi$ given by:
\[\psi_0=\left\{\begin{array}{l}-1:H^{2k}=\Z \to H_{2k+1}=\Z\\1:H^{2k+1}=\Z
\to H_{2k}=\Z\end{array}\right.\;\;\;\psi_1=1:H^{2k+1}\to H_{2k+1}\]
Geometrically $(H,\psi)$ is a symmetric Poincar\'e complex over $\Z$ which is
cobordant to the complex associated to the mapping torus of the
self-diffeomorphism of $\mathbf{CP}^2$ given
by complex conjugation.
Again by direct computation, $\chi_{1/2}(H,\Z_2)=1$, $\chi_{1/2}(H,\Q)=0$
 and $\taunew(H,\psi) = 0$.
By considering the absolute torsion of these two generators we see
that the map $L_{rh}^{4k+1}(\Z) \to K_1(\Z)$ is given by:
\[(C,\phi) \mapsto \chi_{1/2}(C,\Q)\tau(-1)\]
thus proving part two of proposition \ref{prop:poincaresign}.

To prove part 1 of proposition \ref{prop:poincaresign} we use the
following lemma taken from \cite{hkr}:

\begin{lemma}
We have the following relationship between \(\taunew\) and signature
modulo 4 of a \(4k\)-dimensional Poincar\'e complex \((C,\phi)\):
\[\sigma(C) = 2\taunew(C,\phi)+(2k+1)\chi(C) \in \Z_4\]
where the map $2:K_1(\Z)=\Z_2 \to \Z_4$ takes $\tau(-1)$ to $2\in\Z_4$.
\end{lemma}

A simple rearrangement of the formula of the above lemma yields the
first part of proposition \ref{prop:poincaresign}.

\bibliographystyle{amsplain}

\end{document}